\documentclass[12pt]{amsart}
\usepackage{amsfonts}
\usepackage{amsmath,amssymb,amsthm}
\usepackage{amstext}
\usepackage{verbatim}
\input{epsf}
\newtheorem{thm}{Theorem}[section]

\newtheorem{lem}[thm]{Lemma}
\newtheorem{prop}[thm]{Proposition}
\newtheorem{prob}[thm]{Problem}
\newtheorem{bprob}[thm]{Basic problem}
\theoremstyle{definition}

\theoremstyle{remark}
\newtheorem{rem}[thm]{Remark}

\begin{document}

\title[On a generalization of criteria \emph{A} and \emph{D} for congruence of triangles]
{On a generalization of criteria \emph{A} and \emph{D} for congruence of triangles}
\author{Vesselka Mihova and Julia Ninova}%
\address{Faculty of Mathematics and Informatics University of Sofia
e-mail: mihova@fmi.uni-sofia.bg}%
\address{Faculty of Mathematics and Informatics University of Sofia
e-mail: julianinova@hotmail.com}

\subjclass{Primary 51F20, Secondary 51M15}
\keywords{congruence of triangles, compare triangles}%

\maketitle
\thispagestyle{empty}

\vskip 4mm

\begin{abstract} The conditions determining that two
triangles are congruent play a basic role in planimetry. By
comparing not congruent triangles with respect to given sets of
corresponding elements it is important to discover if they have any
common geometric properties characterizing them. The present paper
is devoted to an answer of this question. We give a generalization
of congruence criteria \emph{A} and \emph{D} for  triangles and
apply it to prove some selected geometric problems.
\end{abstract}
\vskip 4mm

\section{Introduction}

There are six essential elements of every triangle - three angles
and three sides. The method of constructing a triangle varies according to the facts
which are known about its sides and angles.

It is important to know what is the minimum knowledge about the sides and angles
which is necessary to construct a particular triangle.

Clearly all triangles constructed in the same way with the same data must be identically
equal, i. e. they must be of exactly the same size and shape and their areas must be the same.

Triangles which are equal in all respects are called \emph{congruent triangles}.

The four sets of minimal conditions for two triangles to be congruent are set out in the following geometric
criteria.
\begin{itemize}
\item[] \emph{\textbf{Criterion A}. Two triangles are congruent if two sides and the included angle of
one triangle are respectively equal to two sides and the included angle of the other.}
\vskip 2mm

\item[] \emph{\textbf{Criterion B}. Two triangles are congruent if two angles and a side of
one triangle are respectively equal to two angles and a side of the other.}
\vskip 2mm

\item[] \emph{\textbf{Criterion C}. Two triangles are congruent if the three sides of
one triangle are respectively equal to the three sides of the other.}
\vskip 2mm

\item[] \emph{\textbf{Criterion D}. Two triangles are congruent if two sides and the angle opposite the
greater side of one triangle are respectively equal to two sides and the angle opposite the
greater side of the other.}
\end{itemize}
\vskip 2mm

It should be noted that  in \emph{criteria} \emph{A} and \emph{D}
the sets of corresponding equal elements are two sides and an angle.

In fact the angle given may be any one of the three angles of the
triangle. The problem to ``\emph{Construct a triangle with two of
its sides $\,a\,$ and $\,b,$ $\,a<b,$ and angle  $\,\alpha\,$
opposite the smaller side}" has not a unique solution. There can be
two triangles each of them satisfying the given conditions.

In the present paper we compare not congruent triangles with respect to given sets of
corresponding elements and answer the question what are the geometric properties characterizing such
couples of triangles.

\vskip 4mm
\section{Theoretical basis of the proposed method for comparing triangles}

In $\,\triangle\, ABC\,$ and $\,\triangle\, A_1B_1C_1\,$ it is convenient to use the notations
$\,AB=c,\,$ $BC=a,\,$ $CA=b;$
$\;A_1B_1=c_1,\,$ $B_1C_1=a_1,\,$ $C_1A_1=b_1.$ Let $\,\theta\,$ and $\,\theta_1\,$ be two
corresponding angles of these triangles.
\vskip 2mm

If for $\triangle\, ABC$ and $\triangle\, A_1B_1C_1$ the relations
$a=a_1,\; b=b_1$ and $\theta=\theta_1$ between the corresponding elements
hold, we consider the following cases.
\vskip 2mm

\begin{itemize}
\item The angle $\theta$ is included between the sides $a$ and $b$, i. e. $\theta= \angle ACB$,
and respectively $\theta_1= \angle A_1C_1B_1$. In this case the triangles are congruent
in view of \emph{Criterion A}.

\vskip 2mm
\item Let $a=b$ and correspondingly $a_1=b_1$, i. e. $\triangle\, ABC$ and $\triangle\, A_1B_1C_1$ are isosceles.
Since $\theta=\theta_1$, the triangles are congruent as a
consequence of \emph{Criterion A}.

\vskip 2mm
\item Let $a > b$,  correspondingly $a_1 > b_1$, and the angle $\theta$ is opposite the
greater side $a$. In this case the triangles are congruent in view
of \emph{Criterion D}.

\vskip 2mm
\item Let $a > b$, correspondingly $a_1 > b_1$, and the angle $\theta$ is opposite the
smaller side $b$. In this case the triangles are either congruent or not.
\begin{itemize}
\item[-] If the triangles are congruent, then the angles opposite
the greater sides are necessarily equal.

It could happen that the sum of the equal angles opposite the
greater sides is two right angles. If so, the triangles are
right-angled.
\vskip 2mm

\item[-] If the triangles are not congruent, then we show that the sum of the angles
opposite the greater sides is always two right angles.
\end{itemize}
\end{itemize}

We prove the following
\begin{lem}
Let $\triangle\, ABC$ and $\triangle\, ABD$ be not congruent triangles having $AB$ as a common side.
Let also $AC=AD$. If $\,\angle ABC=\angle ABD$, then $\;\angle ACB+\angle ADB=180^0$.
\end{lem}

\emph{Proof}. Since $\triangle\, ABC$ and $\triangle\, ABD$ are not congruent, then $AC< AB$ (and hence
$AD< AB$). Let us denote $\angle ACB=\alpha$ and $\angle ADB=\beta$.
\begin{figure}[h t b]
\epsfxsize=8cm
\centerline{\epsfbox{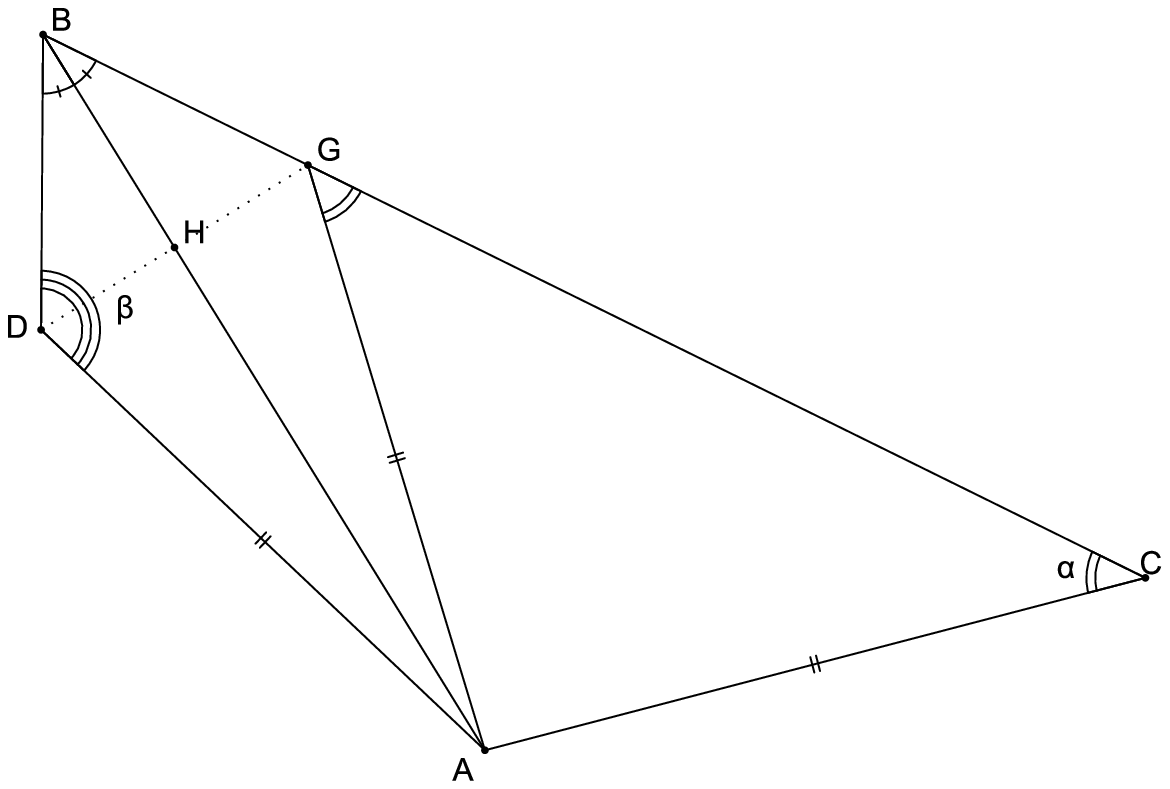}}
\end{figure}
There are two possibilities for the location of $\triangle\, ABC$ and $\triangle\, ABD$ with respect to
the straight line $AB$.
\vskip 2mm

$(i)$ \emph{The points $C$ and $D$ lie on opposite sides of $AB$}.

The symmetry with respect to the straight
line $AB$ transforms $\triangle\, ABD$ into its congruent $\triangle\, ABG$ which lies on one and the same side
of the axis of symmetry $AB$ like $\triangle\, ABC$ (fig. 1). Since  $\triangle\, ABC\ncong \triangle\, ABD$, then
$\triangle\, ABC\ncong \triangle\, ABG$. The condition $\,\angle ABC=\angle ABD$ states that the straight line $AB$
is the bisector of $\angle DBC$. From the symmetry with respect to $AB$ it follows that $G\in BC$  and $BG\neq BC$.
Let, for instance,
$G/BC$. (The case $C/BG$ is analogical.) It is clear that if the conditions of Lemma 2.1 are fulfilled
for $\triangle\, ABC$ and $\triangle\, ABD$, then they are also valid for $\triangle\, ABC$ and $\triangle\, ABG$
and vice versa.

Let us consider $\triangle\, ABC$ and $\triangle\, ABG$. The side $AB$ and $\angle ABC$
are common for both triangles.
In view of the symmetry with respect to $AB$ and $AC=AD$, we get $AD=AG=AC$. Hence, $\triangle\, ACG$ is isosceles
and $\angle ACG=\alpha=\angle AGC$. The angles $\angle AGC$ and $\angle AGB=\angle ADB=\beta$ are adjacent
and hence $\angle AGC+\angle AGB=\angle ACB+\angle ADB=\alpha+\beta=180^0$.

\begin{rem}
The quadrilateral $ACBD$ can be inscribed in a circle.
\end{rem}
\vskip 2mm

$(ii)$ \emph{The points $C$ and $D$ lie on one and the same side of $AB$}.

We consider $\triangle\, ABC$ and $\triangle\, ABG$ (in this case $D\equiv G$). The proof of the
statement is as in $(i)$.
\hfill{$\square$}
\vskip 2mm

\begin{rem}
If $\triangle\, ABC$ and $\triangle\, A_1B_1C_1$ are not congruent, the relations
$AB=A_1B_1,\, AC=A_1C_1$ and $\angle ABC=\angle A_1B_1C_1$ between their corresponding elements
are fulfilled and they have no common side, then
we can choose a suitable congruence and transform $\triangle\, A_1B_1C_1$ into its congruent  $\triangle\, ABD$
so that both triangles  satisfy the conditions of Lemma 2.1.
\end{rem}
\vskip 2mm

Based on the above arguments we can formulate a theorem, which is a generalization of \emph{criteria}
\emph{A} and \emph{D}
for congruence of triangles (see also \cite{Sch}, p. 12).
\vskip 2mm

The denotations $AB=c,\,$ $BC=a,\,$ $CA=b;$
$A_1B_1=c_1,\,$ $B_1C_1=a_1,\,$ $C_1A_1=b_1$ are usually used in $\triangle \,ABC$ and $\triangle\, A_1B_1C_1$
respectively.

\begin{thm}
Let $\theta$ and $\theta_1$ be two corresponding angles of $\triangle \,ABC$ and $\triangle\, A_1B_1C_1$.
If $a=a_1,\, b=b_1$ and $\theta=\theta_1$, then $\triangle\, ABC$ and $\triangle\, A_1B_1C_1$ are either
congruent, or not congruent but the sum of the other two angles, not included between the given sides,
is two right angles.
\end{thm}
\vskip 2mm

Lemma 2.1 and Theorem 2.4 can be used as alternative methods of comparing different triangles.

\vskip 4mm

\section{Application of Theorem 2.4 to two geometric problems}

The solutions of next selected problems are based on Theorem 2.4.

\begin{prob}
$($\cite{NM2}, problems $4.20$ and $4.23$; \cite{T}$)$ Let the middle points of the sides
$BC$, $CA$ and $AB$ of $\,\triangle \,ABC$ be
$F$, $D$, and $E$ respectively. If the center $G$ of the circumscribing
circle $k$ of $\,\triangle\, FDE$ lies on the bisector of $\angle ACB$, then
$\triangle\, ABC$ is either isosceles $(CA=CB)$, or not isosceles but $\angle ACB=60^0$.
\end{prob}

\emph{Proof}. It is given that the center $G$ of the circumscribing
circle $k$ of $\triangle\, FDE$ lies on the bisector of $\angle ACB$ (fig. 2).

\begin{figure}[h t b]
\epsfxsize=8cm
\centerline{\epsfbox{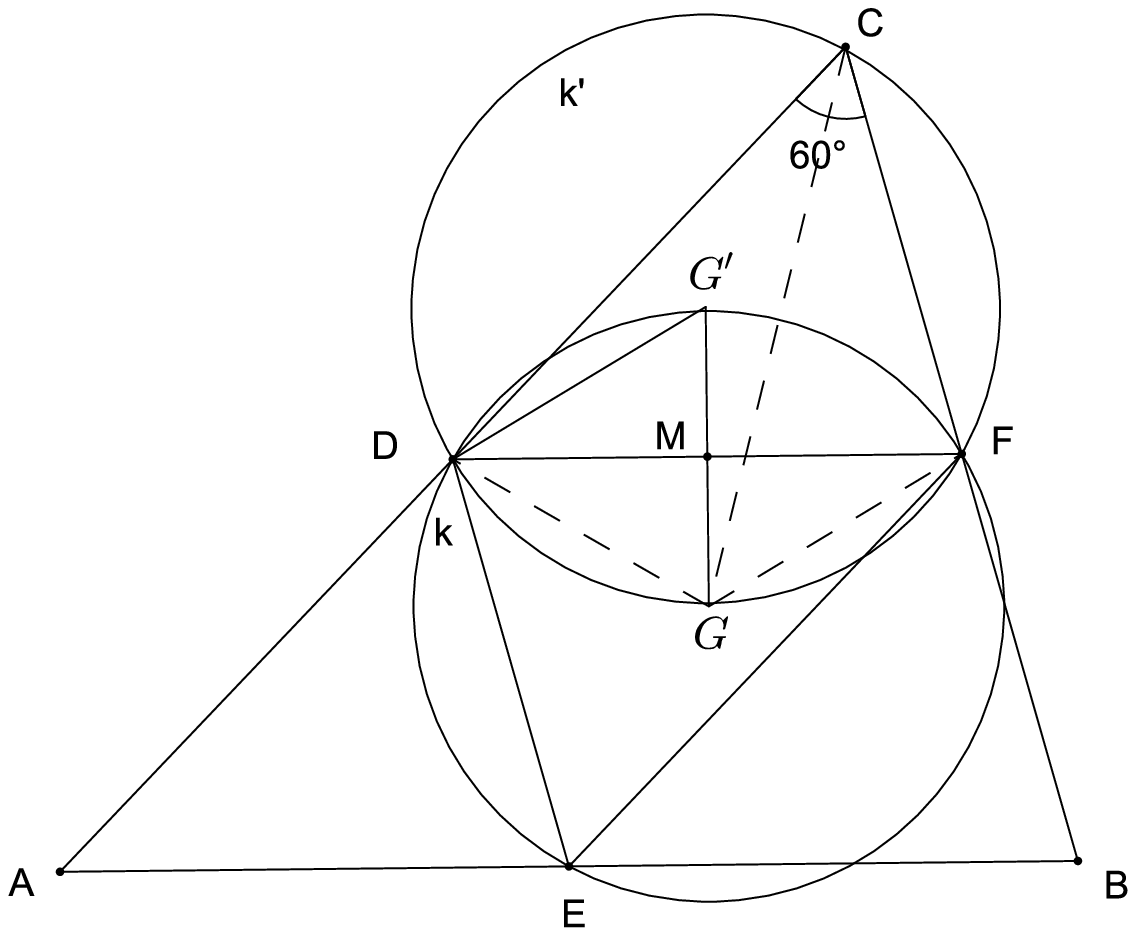}}
\end{figure}
Since $\triangle\, CGD$ and $\triangle \,CGF$ have a common side $CG$, equal corresponding angles
$\angle DCG=\angle FCG$ and equal corresponding sides $DG=FG$ (as radii of $k$), the conditions of
Theorem 2.4 are satisfied. Then $\triangle\, CGD$ and $\triangle\, CGF$ are either congruent, or not congruent.
\vskip 2mm

$(i)$ If $\triangle\, CGD$ and $\triangle\, CGF$ are congruent, then $CD=CF$ and hence $CA=CB$, i. e.
$\triangle \,ABC$ is isosceles.

\begin{rem}
There are two possibilities for $\angle ACB$:
\begin{itemize}
\item[-] either $\angle ACB=60^0$ and $\triangle\, ABC$ is equilateral,
\vskip 2mm

\item[-] or $\angle ACB\neq 60^0$ and $\triangle\, ABC$ is isosceles.
\end{itemize}
\end{rem}
\vskip 2mm

$(ii)$ If $\triangle\, CGD$ and $\triangle\, CGF$ are not congruent, then in view of Lemma 2.1
$\angle CDG+\angle CFG=180^0$ and the quadrilateral $CDGF$ can be inscribed in a circle $k'$ (fig. 2).

It is easy to be seen that $\triangle\, EFD\cong \triangle\, CDF$ and the circumscribing circles $k$ and $k'$
have equal radii. The circles $k$ and $k'$ are symmetrically located with respect to their common chord $FD$.
Since the center $G$ of $k$ lies on $k'$, then the center $G'$ of $k'$ lies on $k$. Hence,
$\triangle\, DGG'\cong\triangle\, FGG'$, both triangles are equilateral, $\angle DGF=120^0$ and $\angle ACB=60^0$.
\hfill{$\square$}

\vskip 4mm
\begin{prob} $($\cite{NM1}, Problem $8$; \cite{NM2}, Problem $4.12)$
Let in $\triangle \,ABC$ the straight lines $AA_1,\,A_1\in BC$,
and $BB_1,\, B_1\in AC$, be the bisectors of $\angle CAB$ and $\angle CBA$ respectively.
Let also  $AA_1\cap BB_1=J$. If $JA_1=JB_1$
then $\triangle \,ABC$ is either isosceles $\,(CA = CB)$, or not isosceles but $\angle ACB = 60^0$.
\end{prob}

\emph{Proof}. We use the denotations $\angle BAC=2\alpha$, $\angle
ABC=2 \beta$, $\angle ACB=2\gamma$.

Since $J$ is the cut point of the angle bisectors $\,AA_1\,$ and $\,BB_1\,$ of $\triangle\, ABC$,
then the straight line $CJ$ is the bisector of $\angle ACB$ and
$\alpha+\beta+\gamma=90^0$ (fig. 3).
\begin{figure}[h t b]
\epsfxsize=8cm
\centerline{\epsfbox{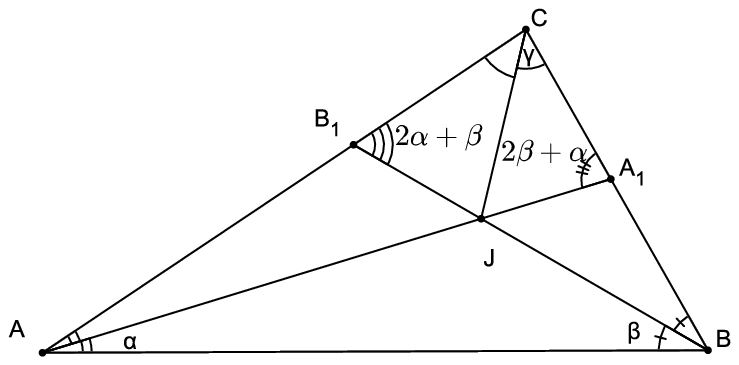}}
\end{figure}

Since $\angle CB_1J$ is an exterior angle of
$\triangle\, ABB_1$, then $\angle CB_1J=2\alpha+\beta$.
Since $\angle CA_1J$ is an exterior angle of
$\triangle\, ABA_1$, then $\angle CA_1J=2\beta+\alpha$.

Let us compare $\triangle\, CA_1J$ and $\triangle\, CB_1J$. They have a common side $CJ$,
corresponding equal sides $JA_1=JB_1$ and angles $\angle A_1CJ=\angle B_1CJ$.

The conditions of Theorem 2.4 are satisfied. Then $\triangle\, CA_1J$ and $\triangle\, CB_1J$
are either congruent, or not.
\vskip 2mm

$(i)$ If $\triangle\, CA_1J$ and $\triangle\, CB_1J$ are congruent, then their corresponding elements
are equal, in particular
$$\angle CB_1J=\angle CA_1J\;\Leftrightarrow\; 2\alpha+\beta=2\beta+\alpha\;\Leftrightarrow\;
\alpha=\beta.$$

Hence, $\triangle\, ABC$ is isosceles with $CA=CB$.
\vskip 2mm

\begin{rem}
There are two possibilities for $\angle ACB$:
\begin{itemize}
\item[-] $\angle ACB=60^0$ and $\triangle\, ABC$ is equilateral;
\vskip 2mm

\item[-] $\angle ACB\neq 60^0$ and $\triangle\, ABC$ is isosceles.
\end{itemize}
\end{rem}
\vskip 2mm

$(ii)$ If $\triangle\, CA_1J$ and $\triangle\, CB_1J$ are not congruent, then with respect to
Lemma 2.1
$$\angle CB_1J+\angle CA_1J=180^0 \;\Leftrightarrow\; (2\alpha+\beta)+(2\beta+\alpha)=180^0 \;\Leftrightarrow\;
\alpha+\beta=60^0.$$

Hence, $\angle ACB=60^0$.
\hfill{$\square$}

\vskip4mm
\section{Groups of problems}
In this section we illustrate the composing technology
of new problems as an interpretation of  specific logical models.

Our aim is the \emph{basic problem} in each of the groups under consideration to be with (exclusive or not
exclusive) disjunction as a logical structure in the conclusion and its proof to be based
on Lemma 2.1 or on Theorem 2.4.

\vskip 4mm

\subsection{Problems of group I}

Suitable logical models for formulation of \emph{equivalent} problems
and \emph{generating} problems of a given problem are described in detail in \cite {NM1, NM2}.
\vskip 2mm

The basic statements we need in this group of problems are:
\begin{itemize}
\item[] $t:=\{$ \it A square with center $O$ is inscribed in a $\triangle\, ABC$ in the following way:
the vertexes of the square lie on the sides of the triangle, in addition two of them lie on the side $AB$.$\}$
\vskip 2mm

\item[] $p:=\{\angle ACB=90^0\}$
\vskip 2mm

\item[] $q:=\{CA=CB\}$
\vskip 2mm

\item[] $r:=\{\angle ACO=\angle BCO\}$
\end{itemize}
\vskip 2mm

We describe the logical scheme for the composition of the Basic problem 4.4, which has not
exclusive disjunction as a logical structure in the conclusion:

\begin{itemize}
\item[-] First we formulate and prove the \emph{generating}  problems - Problem 4.1 with logical structure
$t\wedge p\rightarrow r$ and Problem 4.3 with logical structure $t\wedge q\rightarrow r$.
\vskip 2mm

\item[-] To generate problems with logical structure $\; (*)\quad t\wedge(p \vee q)\rightarrow r$
we use the logical equivalence
$$(t\wedge p\rightarrow r)\wedge (t\wedge q\rightarrow r)\;\Leftrightarrow\; t\wedge (p\vee q) \rightarrow r.$$
\vskip 2mm

\item[-] Finally, the formulated \emph{inverse} problem - Basic problem 4.4 - to the problem with structure $(*)$
has the logical structure $t\wedge r\rightarrow p \vee q$.
\end{itemize}

\begin{prob}
In $\triangle\, ABC$ is inscribed a square with center $O$ in the
following way: the vertexes of the square lie on the sides of the
triangle, in addition two of them lie on the side $AB$. Prove that
if $\angle ACB=90^0$, then $\angle ACO=\angle BCO$.
\end{prob}

\emph{Proof}.  Let the quadrilateral $MNPQ$, $M\in AB$, $\,N\in AB$, $\, P\in BC$, $\,Q\in AC$,
be the inscribed in $\triangle\, ABC$ square (fig. 4).
\begin{figure}[h t b]
\epsfxsize=8cm
\centerline{\epsfbox{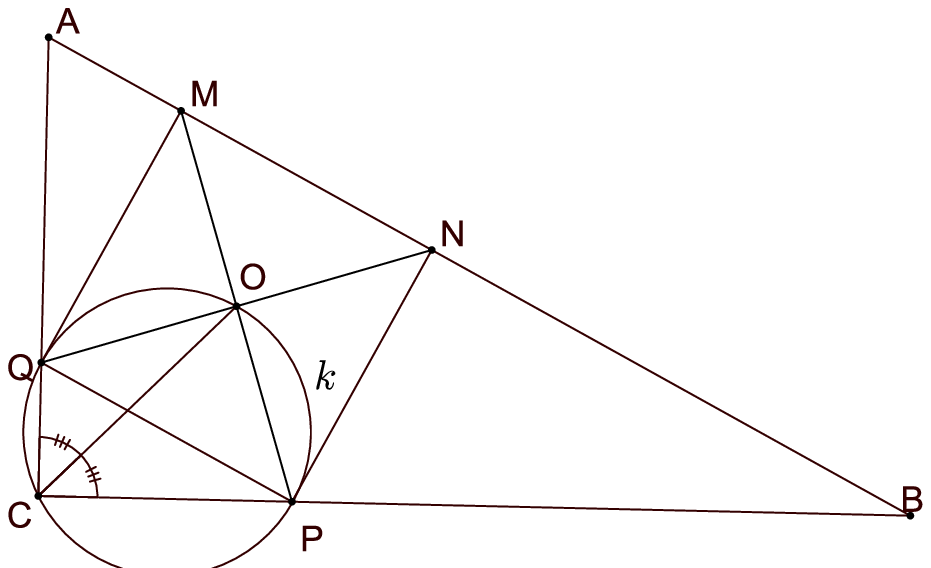}}
\end{figure}
Since the diagonals of a square are equal, intersect at right angles, bisect each other and
bisect the opposite angles, then $OP=OQ$ and $\angle POQ=90^0$. The quadrilateral $OPCQ$
can be inscribed in a circle $k$ with diameter $PQ$. To the equal chords $OQ$ and $OP$ of $k$
correspond equal angles, i. e. $\angle ACO=\angle BCO$.
\hfill{$\square$}

\vskip 3mm
\begin{prob}
In $\triangle\, ABC$ is inscribed a rectangle with center $O$ in the
following way: the vertexes of the rectangle lie on the sides of the
triangle, in addition two of them lie on the side $AB$. Prove that
if $CA=CB$, then $\angle ACO=\angle BCO$.
\end{prob}

\emph{Proof}.  Let the quadrilateral $MNPQ$, $M\in AB$, $\,N\in AB$, $\, P\in BC$, $\,Q\in AC$,
be the inscribed in $\triangle\, ABC$ rectangle (fig. 5).
\begin{figure}[h t b]
\epsfxsize=8cm
\centerline{\epsfbox{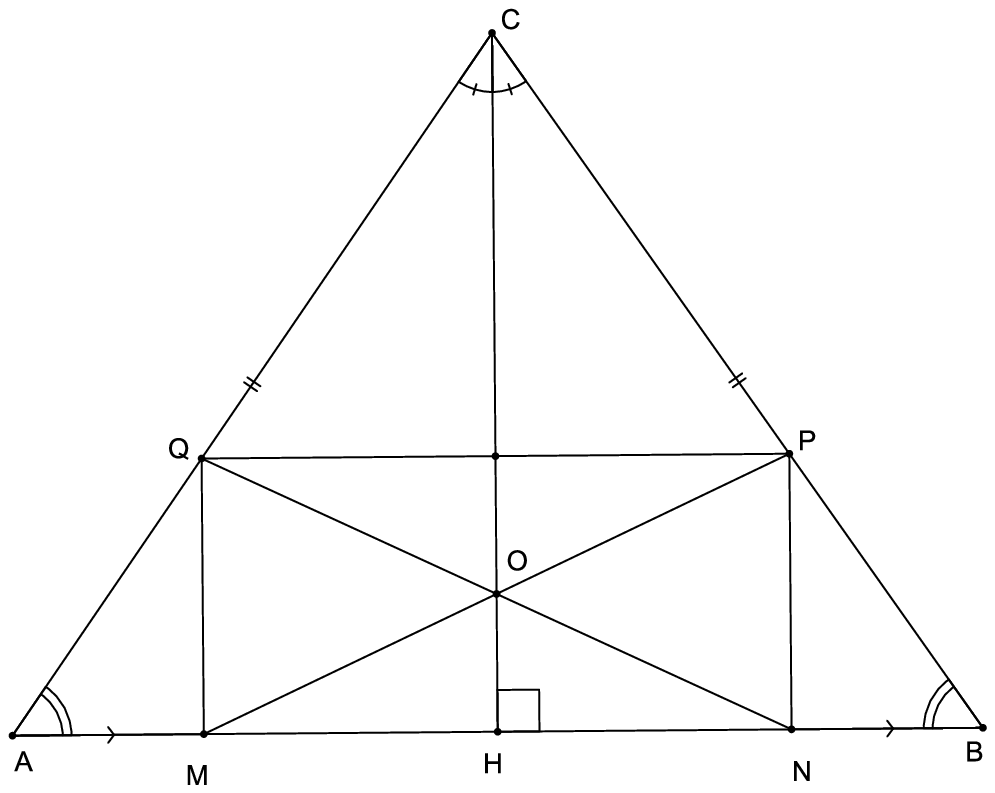}}
\end{figure}
Since the diagonals of a rectangle are equal and bisect each other, then $OM=ON=OP=OQ$.

Let $CH\perp AB,\, H\in AB$. Provided that $\triangle\, ABC$ is isosceles with $CA=CB$, the point
$H$ is the middle point of $AB$ and the straight line $CH$ is the bisector of  $\angle ACB$.

Because $MQ \parallel NP,\; NP \parallel CH$ and $MQ=NP$, it follows that $\triangle\, AMQ \cong \triangle\, BNP$
(by \emph{Criterion B}) and $AM=BN$. Hence, $H$ is also the middle point of $MN$. Since
$\triangle\, MON$ is isosceles, then its median $OH$ is also an altitude, i. e.
$OH\perp MN$. This means that $O\in CH$ and $\angle ACO=\angle BCO$.
\hfill{$\square$}
\vskip 2mm

A special case of Problem 4.2 is  Problem 4.3  with a logical structure $t\wedge q\,\rightarrow \,r$.
\vskip 2mm

\begin{prob}
In $\triangle\, ABC$ is inscribed a square with center $O$ in the
following way: the vertexes of the square lie on the sides of the
triangle, in addition two of them lie on the side $AB$. Prove that
if $CA=CB$, then $\angle ACO=\angle BCO$.
\end{prob}
\vskip 4mm

Now we formulate and prove the \emph{Basic problem} in this group.

\begin{bprob}
In $\triangle\, ABC$ is inscribed a square with center $O$ in the following way:
the vertexes of the square lie on the sides of the triangle, in addition two of them lie on the side $AB$.
Prove that if $\angle ACO=\angle BCO$, then $CA=CB$ or $\angle ACB=90^0$.
\end{bprob}

\emph{Proof}. Let the quadrilateral $MNPQ$, $M\in AB$, $\,N\in AB$, $\, P\in BC$, $\,Q\in AC$,
be the inscribed in $\triangle\, ABC$ square (fig. 6).
\begin{figure}[h t b]
\epsfxsize=8cm
\centerline{\epsfbox{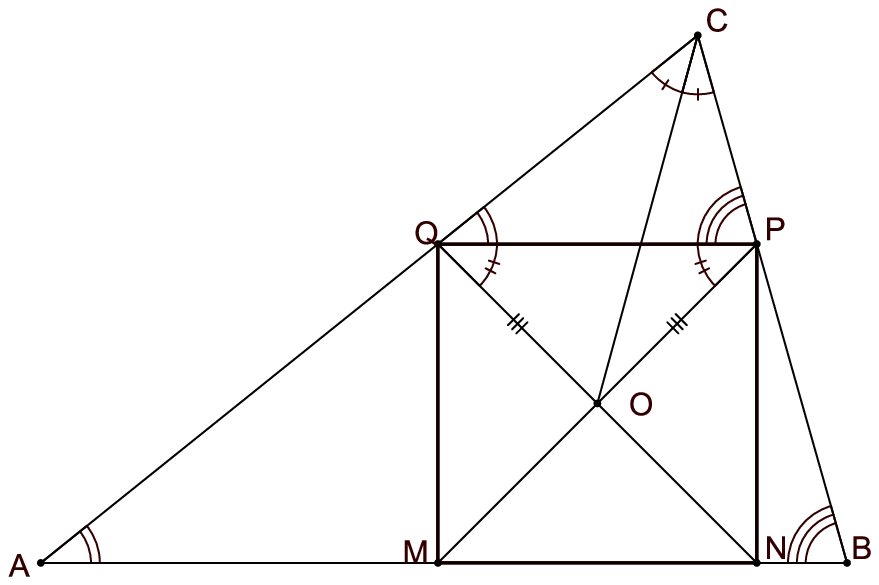}}
\end{figure}
Since the diagonals of any square are equal, intersect at right angles, bisect each other and
bisect the opposite angles, then $OP=OQ$ and $\angle OPQ=\angle OQP=45^0$.

We compare $\triangle\, CQO$ and $\triangle\, CPO$. They have a common side $CO$, respectively equal sides
$OQ=OP$ and angles $\angle QCO=\angle PCO$. We compute that $\angle CQO=\angle CAB+45^0$
and $\angle CPO=\angle CBA+45^0$ as exterior angles of $\triangle\, QAN$ and $\triangle\, PBM$
respectively.

In view of Theorem 2.4 $\;\triangle\, CQO$ and $\triangle\, CPO$ are either congruent or not.

\begin{itemize}
\item[(i)] If $\triangle\, CQO$ and $\triangle\, CPO$ are congruent, then $\angle CQO=\angle CPO$ and hence
$\angle CAB=\angle CBA$, i. e. $CA=CB$ and $\triangle\, ABC$ is isosceles.

In this case $\angle ACB$ is either a right angle and $\triangle\, ABC$ is isosceles
right-angled, or not a right angle and $\triangle\, ABC$ is only isosceles.
\vskip 2mm

\item[(ii)] If $\triangle\, CQO$ and $\triangle\, CPO$ are not congruent then, in view of Lemma 2.1,
$\angle CQO+\angle CPO=180^0$ and hence $\angle CAB+\angle CBA=90^0$, i. e. $\triangle\, ABC$ is
right-angled with $\angle ACB=90^0$.
\end{itemize}
\hfill{$\square$}

\begin{rem}
 A logically incorrect version of the Basic problem 4.4 is Problem 1.54 in \cite{KL}.
\end{rem}
\vskip 2mm

We reformulate Problem 4.4 by keeping the condition of homogeneity of the conclusion.
\vskip 2mm

\begin{prob}
In $\triangle\, ABC$ is inscribed a square with center $O$ in the following way:
the vertexes of the square lie on the sides of the triangle, in addition two of them lie on the side $AB$.
Prove that if $\angle ACO=\angle BCO$, then $\triangle\, ABC$ is either isosceles with $CA=CB$,
or not isosceles but right-angled  with $\angle ACB=90^0$.
\end{prob}
\vskip 4mm

\subsection{Problems of group II}
\vskip 3mm

By formulating appropriate statements and giving suitable logical models we get
two \emph{generating} problems that are necessary for the construction of the Basic problem 4.9.

The basic statements we need are:
\begin{itemize}
\item[] $t:=\{$\it In $\triangle \,ABC$ the straight lines $AA_1,\,A_1\in BC$,
and $BB_1,\, B_1\in AC$, are the bisectors of $\angle CAB$ and $\angle CBA$ respectively.$\}$
\vskip 2mm

\item[] $p:=\{\angle ACB=60^0\}$
\vskip 2mm

\item[] $q:=\{\angle CAB=120^0\}$
\vskip 2mm

\item[] $r:=\{\angle BB_1A_1=30^0\}$
\end{itemize}
\vskip 2mm

Since the sum of the angles of any triangle is equal to two right angles,
the statements $\,p\,$ and $\,q\,$ are mutually exclusive. Hence, if $\,p\,$ is true,
so is  $\,\neg q\,$ and vice versa.
\vskip 4mm

We describe the logical scheme for the composition of the Basic problem 4.9, which has exclusive
disjunction as a logical structure in the conclusion:
\begin{itemize}
\item[-] First we formulate and prove the \emph{generating}  problems - Problem 4.7 with logical structure
$t\wedge p\rightarrow r$ and Problem 4.8 with logical structure $t\wedge q\rightarrow r$.
\vskip 2mm

\item[-] Since the statements $\,p\,$ and $\,q\,$ are mutually exclusive then the equivalences
$\;p\wedge \neg q\;\Leftrightarrow\; p$ and $\; \neg p\wedge q \;\Leftrightarrow\; q\;$ are true.
As a consequence of these facts problems with logical structures $\;t\wedge p\;\rightarrow\; r\;$ and
$\;t\wedge (p\wedge \neg q)\;\rightarrow\; r\;$ are equivalent. So the problems with logical structures
$\;t\wedge q\;\rightarrow\; r\;$ and
$\;t\wedge (q\wedge \neg p)\;\rightarrow\; r$.

To generate problems with logical structure $\; (**)\quad t\wedge(p \veebar q)\rightarrow r$
we use the logical equivalence
$$(t\wedge(p\wedge \neg q)\rightarrow r)\wedge (t\wedge(\neg p\wedge q)\rightarrow r)\quad
\Leftrightarrow\quad t\wedge(p \veebar q)\rightarrow r.$$
\vskip 2mm

\item[-] Finally, the formulated \emph{inverse} problem - the Basic problem 4.9 - to the problem
with structure $(**)$
has the logical structure $t\wedge r\rightarrow p \veebar q$.
\end{itemize}
\vskip 2mm

\begin{prob}
Let in $\triangle \,ABC$ the straight lines $AA_1,\,A_1\in BC$,
and $BB_1,\, B_1\in AC$, be the bisectors of $\angle CAB$ and $\angle CBA$ respectively.
Prove that if $\angle ACB=60^0$, then $\angle BB_1A_1=30^0$.
\end{prob}

\begin{figure}[h t b]
\epsfxsize=8cm
\centerline{\epsfbox{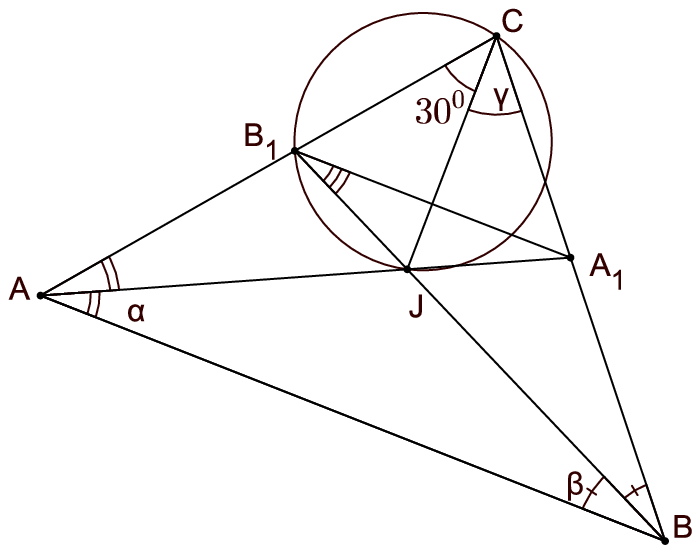}}
\end{figure}

\emph{Proof}. Let $\angle BAA_1=\angle CAA_1=\alpha$, $\,\angle ABB_1=\angle CBB_1=\beta$,
$\, J=AA_1\cap BB_1$.

Since $J$ is the cut point of the angle bisectors of $\triangle\, ABC$,
then $\angle JCA=\angle JCB=\gamma=30^0$ (fig. 7).

Because $\alpha+\beta+\gamma=90^0$
it follows that $\angle AJB=120^0$. Hence, the quadrilateral $CA_1JB_1$ can be inscribed in a circle.
Then  $\angle JA_1B_1=\angle JCB_1=30^0$ and $\;\angle JB_1A_1=\angle JCA_1=30^0$
as angles in the same segment of this circle.
\hfill{$\square$}
\vskip 4mm

\begin{prob}
Let in $\triangle \,ABC$ the straight lines $AA_1,\,A_1\in BC$,
and $BB_1,\, B_1\in AC$, be the bisectors of $\angle CAB$ and $\angle CBA$ respectively.
Prove that if $\angle BAC=120^0$, then $\angle BB_1A_1=30^0$.
\end{prob}

\emph{Proof}. Let $J=AA_1\cap BB_1$, $\, E=A_1B_1\cap CJ$, $\, C_1=CJ\cap AB$.
Since $\angle BAC=120^0$, then its adjacent angles have a measure of $60^0$. It is easy to be seen that
the point $B_1$ is equidistant from the straight lines $BA$, $\,BC$, $\,AA_1$ and that the straight line
$A_1B_1$ is the bisector of $\angle CA_1A$ (fig. 8).
\begin{figure}[h t b]
\epsfxsize=8cm
\centerline{\epsfbox{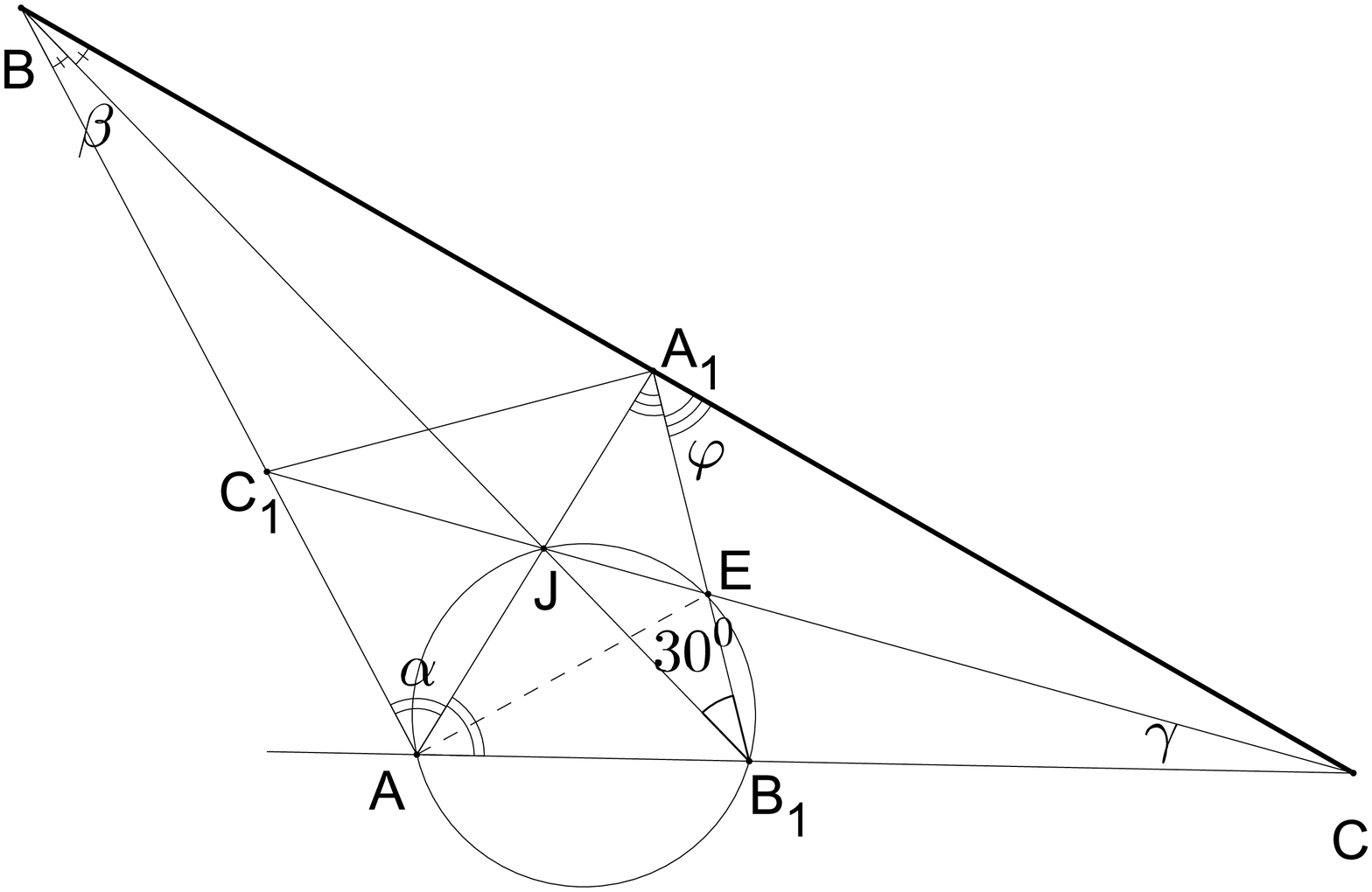}}
\end{figure}
The proof that the straight line $A_1C_1$ is the bisector of $\angle BA_1A$ is analogical.
It follows that $\angle B_1A_1C_1$ is a right angle (the bisectors of any two adjacent angles
are perpendicular to each other) (see also \cite {M}, p. 194, Problem 156).

As a consequence we get that $E$ is the intersection point of the angle bisectors $CJ$ and $A_1B_1$
of $\triangle\, AA_1C$ and hence $\angle JAE=\angle EAB_1=30^0$.

Let $\varphi=\angle CA_1B_1=\angle B_1A_1A$ and $\gamma=\angle C_1CA=\angle C_1CB$.
Then $\angle A_1B_1C=60^0+\varphi$ as an exterior angle of  $\triangle\, A_1B_1A$,
the sum of the angles of $\triangle\, AA_1C$ is $60^0+2\varphi+2\gamma=180^0$, i. e.
$\varphi+\gamma=60^0$ and hence $\angle JEB_1=120^0$.

Thus, the quadrilateral $AJEB_1$ can be inscribed in a circle.
We conclude that $\angle JAE=\angle JB_1E=30^0$ as angles in the same segment of this circle.
Hence, $\angle BB_1A_1=30^0$.
\hfill{$\square$}
\vskip 4mm

Now we formulate and prove the \emph{Basic problem} in this group.

\begin{bprob}
 Let in $\triangle \,ABC$ the straight lines $AA_1,\,A_1\in BC$,
and $BB_1,\, B_1\in AC$, be the bisectors of $\angle CAB$ and $\angle CBA$ respectively.
Prove that if
$\,\angle BB_1A_1= 30^0$, then either $\angle ACB=60^0$ or $\angle BAC=120^0$.
\end{bprob}

\emph{Proof}. Let us denote $\angle BAA_1=\angle CAA_1=\alpha$, $\;\angle ABB_1=\angle CBB_1=\beta$,
 $AA_1\cap BB_1=J$.

 Since $J$ is the cut point of the angle bisectors of $\triangle\, ABC$,
then the straight line $CJ$ is the bisector of $\angle ACB$. Denoting
$\;\gamma=\angle JCA=\angle JCB\;$ we get $\alpha+\beta+\gamma=90^0$ (fig. 9).
\begin{figure}[h t b]
\epsfxsize=8cm
\centerline{\epsfbox{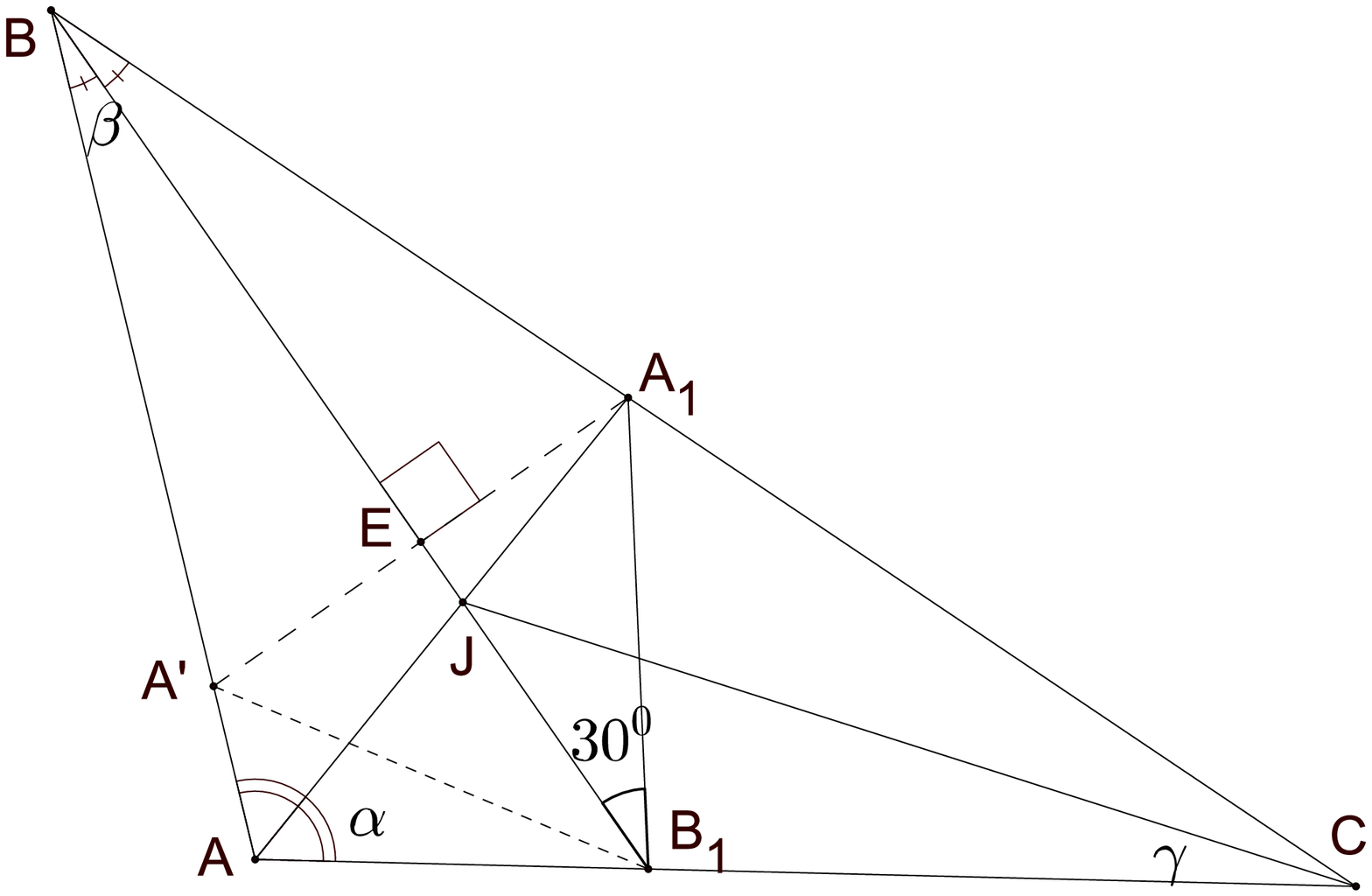}}
\end{figure}

Let the point $A'$ be orthogonally symmetric to the point $A_1$ with respect to the axis $BB_1$.
It follows that $A'\neq A$. (If $A'\equiv A$ then $\triangle\, ABC$ does not exist.)
The straight line $BB_1$ is the bisector of $\angle ABC$ and consequently $A'\in AB$ and $B_1A_1=B_1A'$.
On the other hand, $\angle BB_1A_1= 30^0$ and hence $\triangle\, A_1B_1A'$ is equilateral.

We compute $\angle AA'B_1=30^0+\beta$ (as an exterior angle of $\triangle\, A'BB_1$),
$\,\angle AA'A_1=90^0+\beta$ (as an exterior angle of $\triangle\, A'BE$),
$\,\angle AB_1A' =60^0+\gamma-\alpha\,$ and $\,\angle AB_1A_1=120^0+\gamma-\alpha$.

Let us compare $\triangle\, AA_1B_1$ and $\triangle\, AA_1A'$. They have a common side $AA_1$ and
corresponding equal sides $A_1B_1=A_1A'$ and angles $\angle B_1AA_1=\angle A'AA_1=\alpha$. In view
of Theorem 2.4 we have the possibilities:
\begin{itemize}
\item[(i)] $\triangle\, AA_1B_1$ and $\triangle\, AA_1A'$ are congruent. Then
$\angle AB_1A_1=\angle AA'A_1$, i. e.
$120^0+\gamma-\alpha=90^0+\beta$. Hence, $2 \gamma=\angle ACB=60^0$.
\vskip 2mm

\item[(ii)] $\triangle\, AA_1B_1$ and $\triangle\, AA_1A'$ are not congruent. By Lemma 2.1
it follows that $\angle AB_1A_1+\angle AA'A_1=180^0$, i. e.
$(120^0+\gamma-\alpha)+(90^0+\beta)=180^0$. Hence, $2 \alpha=\angle BAC=120^0$.
\end{itemize}
\hfill{$\square$}
\vskip 2mm

\begin{rem}
An alternate version of Problem 4.9 is  Problem 6, p. 12, in \cite{Sch}.
\end{rem}
\vskip 2mm

In order to formulate a special type equivalent problem to this Basic problem we prove
\begin{prop}
If the statements $p$ and $q$ are mutually exclusive then the following equivalences
are true
$$(\neg(p\veebar q))\;\Leftrightarrow\;
(p\vee \neg q)\wedge (\neg p\vee q)\;\Leftrightarrow\; \neg p\wedge\neg q.$$
\end{prop}
\vskip 2mm

\emph{Proof.}
$$(\neg(p\veebar q))\;\Leftrightarrow\;
\neg((p\wedge\neg q)\vee(\neg p\wedge q))$$

$$\Leftrightarrow\; (p\vee \neg q)\wedge (\neg p\vee q)\;\Leftrightarrow\;
p\wedge (\neg p\vee q)\vee\neg q\wedge (\neg p\vee q)$$

$$\Leftrightarrow\; (p\wedge \neg p)\vee (p\wedge q) \vee (\neg q\wedge\neg p)
\vee (q\wedge \neg q)\;\Leftrightarrow\; \neg p\wedge\neg q.$$
\hfill{$\square$}
\vskip 2mm

Because of this Proposition problems with logical structures
$\; t\wedge (\neg(p\veebar q))\;\rightarrow\; \neg r\;$ and
$\;t\wedge (\neg p\wedge\neg q)\;\rightarrow\; \neg r\;$ are equivalent.

\vskip 2mm

The following problem is equivalent to the Basic problem 4.9.
\begin{prob}
Let in $\triangle \,ABC$ the straight lines $AA_1,\,A_1\in BC$,
and $BB_1,\, B_1\in AC$, be the bisectors of $\angle CAB$ and $\angle CBA$ respectively.
Prove that if $\angle ACB\neq 60^0$ and $\angle CAB\neq 120^0$ then $\angle BB_1A_1\neq 30^0$.
\end{prob}

\emph{Proof}. Assuming the truth of the contrary statement, i. e. $\angle BB_1A_1 = 30^0$,
the solution of this problem leads to the solution of the Basic problem 4.9.
\hfill{$\square$}

\vskip 6mm
Acknowledgements. The first author is partially supported by Sofia University Grant 99/2013.
The second author is partially supported by Sofia University Grant 159/2013.

\end{document}